\documentclass[10pt,a4paper]{article}
\usepackage{amssymb}
\usepackage{amscd}
\usepackage{latexsym}
\usepackage{amsmath}
\usepackage{amsfonts}
\usepackage{amscd}

\begin{document}
\newtheorem{theorem}{Theorem}[section]
\newtheorem{remark}[theorem]{Remark}
\newtheorem{mtheorem}[theorem]{Main Theorem}
\newtheorem{bbtheo}[theorem]{The Strong Black Box}
\newtheorem{observation}[theorem]{Observation}
\newtheorem{proposition}[theorem]{Proposition}
\newtheorem{lemma}[theorem]{Lemma}
\newtheorem{testlemma}[theorem]{Test Lemma}
\newtheorem{mlemma}[theorem]{Main Lemma}
\newtheorem{note}[theorem]{{\bf Note}}
\newtheorem{steplemma}[theorem]{Step Lemma}
\newtheorem{corollary}[theorem]{Corollary}
\newtheorem{notation}[theorem]{Notation}
\newtheorem{example}[theorem]{Example}
\newtheorem{definition}[theorem]{Definition}

\renewcommand{\labelenumi}{(\roman{enumi})}
\def\Pf{\smallskip\goodbreak{\sl Proof. }}

\def\Fin{\mathop{\rm Fin}\nolimits}
\def\br{\mathop{\rm br}\nolimits}
\def\fin{\mathop{\rm fin}\nolimits}
\def\Ann{\mathop{\rm Ann}\nolimits}
\def\Aut{\mathop{\rm Aut}\nolimits}
\def\End{\mathop{\rm End}\nolimits}
\def\bfb{\mathop{\rm\bf b}\nolimits}
\def\bfi{\mathop{\rm\bf i}\nolimits}
\def\bfj{\mathop{\rm\bf j}\nolimits}
\def\df{{\rm df}}
\def\bfk{\mathop{\rm\bf k}\nolimits}
\def\bEnd{\mathop{\rm\bf End}\nolimits}
\def\iso{\mathop{\rm Iso}\nolimits}
\def\id{\mathop{\rm id}\nolimits}
\def\Ext{\mathop{\rm Ext}\nolimits}
\def\Ines{\mathop{\rm Ines}\nolimits}
\def\Hom{\mathop{\rm Hom}\nolimits}
\def\bHom{\mathop{\rm\bf Hom}\nolimits}
\def\Rk{ R_\k-\mathop{\bf Mod}}
\def\Rn{ R_n-\mathop{\bf Mod}}
\def\map{\mathop{\rm map}\nolimits}
\def\cf{\mathop{\rm cf}\nolimits}
\def\top{\mathop{\rm top}\nolimits}
\def\Ker{\mathop{\rm Ker}\nolimits}
\def\Bext{\mathop{\rm Bext}\nolimits}
\def\Br{\mathop{\rm Br}\nolimits}
\def\dom{\mathop{\rm Dom}\nolimits}
\def\min{\mathop{\rm min}\nolimits}
\def\im{\mathop{\rm Im}\nolimits}
\def\max{\mathop{\rm max}\nolimits}
\def\rk{\mathop{\rm rk}}
\def\Diam{\diamondsuit}
\def\Z{{\mathbb Z}}
\def\Q{{\mathbb Q}}
\def\N{{\mathbb N}}
\def\bQ{{\bf Q}}
\def\bF{{\bf F}}
\def\bX{{\bf X}}
\def\bY{{\bf Y}}
\def\bHom{{\bf Hom}}
\def\bEnd{{\bf End}}
\def\bS{{\mathbb S}}
\def\AA{{\cal A}}
\def\BB{{\cal B}}
\def\CC{{\cal C}}
\def\DD{{\cal D}}
\def\TT{{\cal T}}
\def\FF{{\cal F}}
\def\GG{{\cal G}}
\def\PP{{\cal P}}
\def\SS{{\cal S}}
\def\XX{{\cal X}}
\def\YY{{\cal Y}}
\def\fS{{\mathfrak S}}
\def\fH{{\mathfrak H}}
\def\fU{{\mathfrak U}}
\def\fW{{\mathfrak W}}
\def\fK{{\mathfrak K}}
\def\PT{{\mathfrak{PT}}}
\def\T{{\mathfrak{T}}}
\def\fX{{\mathfrak X}}
\def\fP{{\mathfrak P}}
\def\X{{\mathfrak X}}
\def\Y{{\mathfrak Y}}
\def\F{{\mathfrak F}}
\def\C{{\mathfrak C}}
\def\B{{\mathfrak B}}
\def\J{{\mathfrak J}}
\def\fN{{\mathfrak N}}
\def\fM{{\mathfrak M}}
\def\Fk{{\F_\k}}
\def\bar{\overline }
\def\Bbar{\bar B}
\def\Cbar{\bar C}
\def\Pbar{\bar P}
\def\etabar{\bar \eta}
\def\Tbar{\bar T}
\def\fbar{\bar f}
\def\nubar{\bar \nu}
\def\rhobar{\bar \rho}
\def\Abar{\bar A}
\def\a{\alpha}
\def\b{\beta}
\def\g{\gamma}
\def\w{\omega}
\def\e{\varepsilon}
\def\o{\omega}
\def\va{\varphi}
\def\k{\kappa}
\def\m{\mu}
\def\n{\nu}
\def\r{\rho}
\def\f{\phi}
\def\hv{\widehat\v}
\def\hF{\widehat F}
\def\v{\varphi}
\def\s{\sigma}
\def\l{\lambda}
\def\lo{\lambda^{\aln}}
\def\d{\delta}
\def\z{\zeta}
\def\th{\theta}
\def\t{\tau}
\def\ale{\aleph_1}
\def\aln{\aleph_0}
\def\Cont{2^{\aln}}
\def\nld{{}^{ n \downarrow }\l}
\def\n+1d{{}^{ n+1 \downarrow }\l}
\def\hsupp#1{[[\,#1\,]]}
\def\size#1{\left|\,#1\,\right|}
\def\Binfhat{\widehat {B_{\infty}}}
\def\Zhat{\widehat \Z}
\def\Mhat{\widehat M}
\def\Rhat{\widehat R}
\def\Phat{\widehat P}
\def\Fhat{\widehat F}
\def\fhat{\widehat f}
\def\Ahat{\widehat A}
\def\Chat{\widehat C}
\def\Ghat{\widehat G}
\def\Bhat{\widehat B}
\def\Btilde{\widetilde B}
\def\Ftilde{\widetilde F}
\def\restr{\mathop{\upharpoonright}}
\def\to{\rightarrow}
\def\arr{\longrightarrow}
\def\LA{\langle}
\def\RA{\rangle}
\newcommand{\norm}[1]{\text{$\parallel\! #1 \!\parallel$}}
\newcommand{\supp}[1]{\text{$\left[ \, #1\, \right]$}}
\def\set#1{\left\{\,#1\,\right\}}
\newcommand{\mb}{\mathbf}
\newcommand{\wt}{\widetilde}
\newcommand{\card}[1]{\mbox{$\left| #1 \right|$}}
\newcommand{\union}{\bigcup}%\limits}
\newcommand{\inters}{\bigcap}%\limits}
\newcommand{\ER}{{\rm E}}
\def\Proof{{\sl Proof.}\quad}
\def\fine{\ \black\vskip.4truecm}
\def\black{\ {\hbox{\vrule width 4pt height 4pt depth
0pt}}}
\def\fine{\ \black\vskip.4truecm}
\long\def\alert#1{\smallskip\line{\hskip\parindent\vrule%
\vbox{\advance\hsize-2\parindent\hrule\smallskip\parindent.4\parindent%
\narrower\noindent#1\smallskip\hrule}\vrule\hfill}\smallskip}

\title{On Abelian Groups Having All Proper \\ Characteristic Subgroups Isomorphic}
\author{Andrey R. Chekhlov \\Faculty of Mathematics and Mechanics, Section of Algebra, \\Tomsk State University, Tomsk 634050, Russia\\{\small e-mails: a.r.che@yandex.ru, cheklov@math.tsu.ru}
\\and\\ Peter V. Danchev \\Institute of Mathematics and Informatics, Section of Algebra, \\Bulgarian Academy of Sciences, Sofia 1113, Bulgaria\\{\small e-mails: danchev@math.bas.bg, pvdanchev@yahoo.com}}
\maketitle

\begin{abstract}{We consider two variants of those Abelian groups with all proper characteristic subgroups isomorphic and give an in-depth study of their basic and specific properties in either parallel or contrast to the Abelian groups with all proper fully invariant subgroups isomorphic, which are studied in details by the current authors in Commun. Algebra (2015). In addition, we also examine those Abelian groups having at least one proper characteristic subgroup isomorphic to the whole group. The established by us results somewhat extend those obtained by Grinshpon-Nikolskaya in Tomsk State Univ. J. Math. \& Mech. (2011, 2012) and in Commun. Algebra (2011), respectively.}
\end{abstract}

\section{Introduction and Definitions}

Throughout this article, let all groups into consideration be {\it additively} written and {\it abelian}. Our notations and terminology from group theory are mainly standard and follow those from \cite{F}, \cite{F1} and \cite{Kap}, respectively. Another useful source on the explored subject is \cite{P} as well. For instance, if $p$ is a prime integer and $G$ is an arbitrary group, $p^nG=\{p^ng~|~ g\in G\}$ denotes the {\it $p^n$-th power subgroup} of $G$ consisting of all elements of $p$-height greater than or equal to $n\in \mathbb{N}$, $G[p^n]=\{g\in G~|~ p^ng=0, n\in \mathbb{N}\}$ denotes the {\it $p^n$-socle} of $G$, and $G_p=\cup_{n<\o} G[p^n]$ denotes the {\it $p$-component} of the {\it torsion part} $tG=\oplus_p G_p$ of $G$.

On the other hand, if $G$ is a torsion-free group and $a\in G$, then let $\chi_G(a)$ denote the {\it characteristic} and let $\t_G(a)$ denote the {\it type} of $a$, respectively. Specifically, the class of equivalence in the set of all characteristics is just called {\it type} and we write $\t$. If $\chi_G(a) \in \t$, then we write $\t_G(a)=\t$, and so $\t(G)=\{\t_G(a) ~ | ~ 0\neq a\in G\}$ is the set of types of all non-zero elements of $G$. The set $G(\t)=\{g\in G ~|~ \t(g)\geqslant \t\}$ forms a pure fully invariant subgroup of the torsion-free group $G$. Recall that a torsion-free group $G$ is called {\it homogeneous} if all its non-zero elements have the same type.

Concerning ring theory, suppose that all rings which we consider are {\it associative} with {\it identity} element. For any ring $R$, the letter $R^{+}$ will denote its {\it additive group}. To simplify the notation and to avoid a risk of confusion, we shall write $\ER (G)$ for the endomorphism ring of $G$ and $\End (G)=\ER (G)^{+}$ for the endomorphism group of $G$.

As usual, a subgroup $F$ of a group $G$ is called {\it fully invariant} if $\phi(F)\subseteq F$ for any $\phi\in \ER (G)$, while if $\phi$ is an invertible endomorphism (= an automorphism), then $F$ is called a {\it characteristic} subgroup.

Classical examples of important fully invariant subgroups of an arbitrary group $G$ are the defined above subgroups $p^nG$ and $G[p^n]$ for any natural $n$ as well as $tG$ and the maximal divisible subgroup $dG$ of $G$; actually $dG$ is a fully invariant direct summand of $G$ (see, for instance, \cite{F}).

To avoid any confusion and misunderstanding, we shall say that a group $G$ has only {\it trivial fully invariant subgroups} if $\{0\}$ and $G$ are the only ones. Same appears for the characteristic subgroups, respectively.

\medskip

The following notions were stated in \cite{CD}.

\medskip

\noindent{\bf Definition 1}. A non-zero group $G$ is said to be an {\it IFI-group} if either it has only trivial fully invariant subgroups, or all its non-trivial fully invariant subgroups are isomorphic otherwise.

\medskip

\noindent{\bf Definition 2}. A non-zero group $G$ is said to be an {\it IC-group} if either it has only trivial characteristic subgroups, or all its non-trivial characteristic subgroups are isomorphic otherwise.

\medskip

Note that Definition 2 implies Definition 1. In other words, any IC-group is an IFI-group; in fact every fully invariant subgroup is characteristic.

\medskip

\noindent{\bf Definition 3}. A non-zero group $G$ is called a {\it strongly IFI-group} if either it has only trivial fully invariant subgroups, or all its non-zero fully invariant subgroups are isomorphic otherwise.

\medskip

\noindent{\bf Definition 4}. A non-zero group $G$ is called a {\it strongly IC-group} if either it has only trivial characteristic subgroups, or all its non-zero characteristic subgroups are isomorphic otherwise.

\medskip

Notice that Definition 4 implies Definition 3.

On another vein, Definition 4 obviously yields Definition 2, but the converse is false. In fact, in \cite{Kap} was constructed a single non-trivial characteristic subgroup of a $2$-group that are pairwise non-isomorphic, thus giving an example of an IC-group which is surely {\it not} a strongly IC-group; in other words, Definition 4 properly implies Definition 2.

\medskip

Our goal in the present paper is to explore some fundamental and exotic properties of the defined above classes of groups, especially the IC-groups and the strongly IC-groups. In addition, we shall investigate even something more, namely the existence of a non-trivial characteristic subgroup of a given group which subgroup is isomorphic to the whole group, calling these groups {\it weakly IC-groups}. Our motivation to do that is to exhibit and compare the characteristically similarities and discrepancies of these group classes.

The major results established by us are formulated and proved in the next section.

\section{Main Theorems and Examples}

For completeness of the exposition and for the reader's convenience, we first and foremost will give a brief retrospection of the most principal results achieved in \cite{CD} concerning IFI-groups and strongly IFI-groups.

As usual, the symbol $\oplus_{m} G=G^{({m})}$ will denote the {\it external} direct sum of $m$ copies of the group $G$, where $m$ is some ordinal (finite or infinite).

\begin{theorem}
Let $G$ be a $p$-group and let $m\geqslant 2$ be an ordinal. Then $G^{({m})}$ is an IFI-group if and only if $G$ is an IC-group if and only if $G$ is an IPI-group.
\end{theorem}

\begin{proposition}\label{tf} Let $G$ be a torsion-free group. Then $G$ is an IFI-group if and only if $G$ is a strongly IFI-group.
\end{proposition}

\begin{lemma}\label{need} (a) A fully invariant subgroup of an IFI-group is an IFI-group.

(b) A fully invariant subgroup of a strongly IFI-group is a strongly IFI-group.

\end{lemma}

\begin{proposition} A non-zero IFI-group is either divisible or reduced.
\end{proposition}

\begin{theorem}\label{main} The following two points hold:

(i) A non-zero group $G$ is an IFI-group if and only if one of the following holds:

\medskip

$\bullet$ For some prime $p$ either $pG=\{0\}$, or $p^2G=\{0\}$ with $r(G)=r(pG)$.

\medskip

$\bullet$ $G$ is a homogeneous torsion-free IFI-group of an idempotent type.

\medskip

(ii) A non-zero torsion group $G$ is a strongly IFI-group if and only if it is an elementary $p$-group for some prime $p$.

\end{theorem}

\begin{proposition}
Every homogeneous fully transitive torsion-free group of an idempotent type is an IFI-group.
\end{proposition}

\begin{corollary} A direct summand of a fully transitive torsion-free IFI-group is again a fully transitive IFI-group.
\end{corollary}

We next continue with the statements and proofs of the chief results. Before doing that, we record the following two hopefully useful assertion.

\begin{lemma} If $G$ is a decomposable fully transitive homogeneous torsion-free group, then all its characteristic subgroups are fully invariant.
\end{lemma}

\Pf Let us decompose $G=A\oplus B$, where $A,B\neq \{0\}$, $H\leq G$ is a characteristic subgroup and $\pi: G\to A$
is the corresponding projection. If $x\in\pi(H)$, $f\in\mathrm{E}(A)$ and $y=f(x)$, then there exists an element $z\in B$ with $\chi(z)=\chi(y)$. So, $\varphi(y)=z$, $\psi(z)=y$ for some $\varphi\in\mathrm{Hom}(A,B)$,
$\psi\in\mathrm{Hom}(B,A)$. Applying Lemma~\ref{01} appointed below, we get that $y\in H$. Employing once again Lemma~\ref{01} (1), we can deduce the wanted full invariability of $H$.
\fine

\begin{lemma}\label{05} If $G$ is a homogeneous torsion-free transitive group, then all its characteristic subgroups are fully invariant.
\end{lemma}

\Pf Choose $0\neq x\in H$ and $f\in\mathrm{E}(G)$, where $H$ is characteristic in $G$. We have $\chi(x)\leq\chi(f(x))$ and $t(x) = t(f(x))$ whenever $f(x)\neq 0$, so that $\chi(kx) = \chi(f(x))$ for some natural number
$k$, whence $\varphi(kx) = f(x)$, where $\varphi\in\mathrm{Aut} (G)$, and so $f(x)\in H$, as needed.
\fine

\subsection{IC-groups}

We start our work here with some elementary but useful observations:

\medskip

$\bullet$ Since the divisible part of any group is necessarily fully invariant, a non-zero IC-group is either divisible or reduced, and hence it follows from Theorem~\ref{main} that it is either a $p$-group or a homogeneous torsion-free group of an idempotent type.

\medskip

$\bullet$ If $G$ is a torsion-free group, then $G$ is an IC-group if, and only if, $G$ is a strongly IC-group (compare with the next subsection). The proof is as Proposition~2.2 from \cite{CD}.

\medskip

$\bullet$ Every homogeneous transitive torsion-free group of an idempotent type is an IC-group -- in fact, in view of the condition being homogeneous, we know that such a group is necessarily fully transitive.

\medskip

$\bullet$ Any divisible group is an IC-group if, and only if, it is a torsion-free group. (This follows automatically from \cite[Proposition 2.7]{CD}.)

\medskip

We now continue with the following technicalities.

\begin{lemma}
(1) If $p^2G=\{0\}$ for some prime $p$, then every characteristic subgroup $H$ in $G$ is fully invariant.

(2) Every torsion IFI-group is an IC-group.
\end{lemma}

\Pf (1) We have $G=G_1\oplus G_2$, where $G_1$ is elementary and $G_2\cong\bigoplus\mathbb{Z}(p^2)$ whenever $pG\neq \{0\}$. If $pH\neq \{0\}$, then as in Proposition~\ref{03} $H=G$. So, let us assume that $pH=\{0\}$. Supposing $\pi: G\to G_1$ is the corresponding projection and $\pi(H)\neq \{0\}$, it follows directly from Lemma~\ref{01} that $H=G[p]$. However, if $\pi(H)=\{0\}$, then applying once again the same Lemma~\ref{01}, we extract that $H=G_2[p]=pG[p]$. So, in both cases, the subgroup $H$ is fully transitive, as asserted.

(2) follows from (1).
\fine

It is an evident fact that if any element of the endomorphism ring $\mathrm{E}(G)$ is an integer multiplied by invertible, then all fully invariant subgroups of a given group are always characteristic.

Furthermore, by usage of \cite[Theorem 2.12]{CD}, if $G$ is a torsion-free IFI-group whose non-zero endomorphisms are monomorphisms, then $G$ is an E-group, so if the rank of the additive group $\mathrm{E}(G)^+$
is finite for such a group $G$, then the tensor product $\mathrm{E}(G)\otimes_{\mathbb Q} \mathbb{Q}$ is a field, and every endomorphism of $G$ is an integer multiplied by invertible (see the paragraph before
Theorem 4.7 in \cite{KMT}). Thus, one can state the following direct consequence:

\begin{corollary}
If $G$ is a torsion-free IFI-group of finite rank whose non-zero endomorphisms are monomorphisms, then $G$ is an IC-group.
\end{corollary}

However, in the case of torsion-free groups of infinite rank the situation is totally different as the next construction shows. Before listing it, we need the following technicalities.

If the additive group $R^+$ of a ring $R$ is a torsion-free group, then we shall say that the ring $R$ is a {\it torsion-free ring}. Likewise, we put $\pi(R)=\{p\in\mathbb{P}\,|\,pR\neq R\}$, where $h_p^G(x)$ is the
$p$-height of an element $x$ in $G$.

\begin{lemma}\label{properties}\cite[Lemma~44.6]{KMT} The following properties of a torsion-free ring $R$ are equivalent:

\medskip

(1) For every $p\in\pi(R)$ and all $a,b\in R$, if $p\,|\,ab$, then (either) $p\,|\,a$ or $p\,|\,b$;

(2) The equality $h_p(ab)=h_p(a)+h_p(b)$ holds for each $p\in\pi(R)$ and all $a,b\in R$;

(3) The ring $R/pR$ has no zero divisors.
\end{lemma}

We are now ready to give our example.

\begin{example}\label{padic} There exists an IFI-group which is {\it not} an IC-group.
\end{example}

\Pf Let $\sigma\in\hat{\mathbb{Z}}_p$ is invertible, but transcendental on $\mathbb{Z}$, and $R$ is a pure subring
in $\hat{\mathbb{Z}}_p$ generated by $1$ and $\sigma$; standardly, $\hat{\mathbb{Z}}_p$ is the ring of $p$-adic integers. If we set $G=R^+$, then it is obviously an $E$-group. We now intend to illustrate that $R$ is a ring of principal ideals, whence $G$ will be an IFI-group by consulting with \cite[Theorem~2.12]{CD}. To that objective, given $I\lhd R$ is an arbitrary ideal. If $\alpha\in R$, then
$p^n\alpha\in\mathbb{Q}_p[\sigma]$, where as usual $\mathbb{Q}_p$ designates the ring of all rational numbers whose denominators are not divided by $p$, and $\mathbb{Q}_p[\sigma]$ is the ring of polynomials of $\sigma$ with coefficients from $\mathbb{Q}_p$. Letting $\gamma$ be the polynomial "of least degree" coming $I$, then one sees that $\gamma R\leq I$. Notice also that $I$ is a $q$-pure ring (without identity) in $R$ for every prime number $q\neq p$. Thus, for every $\beta\in I$, it follows that $p^m\beta\in\gamma R$ for some $m\in\mathbb{N}$. Indeed, if
$\beta=b_0\sigma^k+b_1\sigma^{k-1}+\dots+b_k$, $\gamma=c_0\sigma^s+c_1\sigma^{s-1}+\dots+c_s$, where
$b_i,c_j\in\mathbb{Q}_p$, $k,s\in\mathbb{N}$, $c_0=p^td_0$, where $d_0\in\mathbb{Q}_p\setminus p\mathbb{Q}_p$,
then $k\geq s$ and $p^t\beta-(b_0d_0^{-1}\sigma^{k-s})\gamma$ has degree on $\sigma$ smaller than $k$, so that we can use an induction. Suppose $\gamma_0=p^s\gamma$, where $\gamma_0\in\hat{\mathbb{Z}}_p\setminus p\hat{\mathbb{Z}}_p$.
Hence $\gamma_0 R$ is a pure subring in $R$. In fact, as $R$ is pure in $\hat{\mathbb{Z}}_p$, we have that
$R/pR\cong\hat{\mathbb{Z}}_p/p\hat{\mathbb{Z}}_p$, an application of Lemma~\ref{properties} gives that
$h_p^R(\gamma_0r)=h_p^R(\gamma_0)+h_p^R(r)=h_p^R(r)$, because $h_p^R(\gamma_0)=0$.

Furthermore, it follows from the noted above properties of elements of $I$ that $I\leq \gamma_0R$. But
$\gamma_0R/\gamma R\cong\mathbb{Z}_{p^s}$ and since $I/\gamma R\leq \gamma_0R/\gamma R$, one verifies that $I/\gamma R\cong \mathbb{Z}_{p^l}$ for some $l\leq s$. In particular, $I=(p^l\gamma)R$, i.e., $I$ is a proper principal ideal of $R$, as expected, and hence $G$ is an $IFI$-group, as indicated.

Moreover, note that $\mathrm{Aut}(G)$ is isomorphic to the group $U(\mathbb{Q}_p)$ consisting of all invertible elements of the ring $\mathbb{Q}_p$. To show that, assume $\alpha\in R$ is invertible, and since
$p^n\alpha\in\mathbb{Q}_p[\sigma]$ it must be that $p^n\alpha=a_0\sigma^k+a_1\sigma^{k-1}+\dots+a_k$ for some
$k\in\mathbb{N}$ and $a_0,a_1,\dots,a_k\in\mathbb{Q}_p$. Similarly, $p^m\alpha^{-1}=b_0\sigma^s+b_1\sigma^{s-1}+\dots+b_s$. Therefore,
$(a_0\sigma^k+a_1\sigma^{k-1}+\dots+a_k)(b_0\sigma^s+b_1\sigma^{s-1}+\dots+b_s)=p^{n+m}$, which is actually the desired contradiction that $\sigma$ is transcendental provided $k\neq 0$. Consequently, in $G$ there will exist a characteristic subgroup of rank $1$, but since $G$ is supposed to be infinite rank it definitely will not be an IC-group, as claimed.
\fine

In view of \cite[Corrolary~2.22]{CD}, one can say that if $G$ is an IFI-group, then the cartesian direct sum $G^{(m)}$ is also an IFI-group for any ordinal $m$. However, the following is true:

\begin{example}\label{direct1} There exists such a group $A$ that $G=A^{(m)}$ is an $IC$-group for every ordinal $m>1$, but $A$ is not an $IC$-group. In particular, the direct summand of an $IC$-group need not be an $IC$-group.
\end{example}

\Pf Take the group $A$ from Example~\ref{padic}. Thus, $G$ is an $IFI$-group, because $A$ is an $IFI$-group. Now, suppose that $H\leq G$ is a characteristic subgroup. Since every endomorphism of $G$ is a sum of automorphisms, it must be that $H$ is fully invariant. Consequently, $H\cong G$, i.e., $G$ is an $IC$-group, as asserted.
\fine

\subsection{Strongly IC-groups}

We start our work here with some elementary but helpful observations:

\medskip

$\bullet$ A non-zero torsion group $G$ is a strongly IC-group if, and only if, it is an elementary $p$-group for some prime $p$. (This follows immediately from Theorem~\ref{main} quoted above.)

\medskip

In virtue of Theorem~\ref{main}, one infers that a strongly IFI-group is exactly a strongly IC-group, as well as a torsion-free strongly IC-group is precisely a torsion free IC-group.

\subsection{Weakly IC-groups}

Before starting our work, it is worthwhile noticing that in \cite{GN1}, \cite{GN2} and \cite{GN3} were investigated those Abelian groups (both torsion and torsion-free) having isomorphic proper fully invariant subgroup. In this subsection, we will initiate an examination of the case when we have a proper characteristic subgroup isomorphic to the whole group. For simplicity of the exposition, we shall call these groups just {\it weakly IC-groups}.

\medskip

We now continue with our principal results. Specifically, we proceed by proving with the following statements.

\begin{lemma}\label{01} Let $G=A\oplus B$ and $\pi: G\to A$ be the corresponding projection. If $H$ is a characteristic subgroup in $G$, then:

(1) The inclusion $f\pi(H)\leq H$ holds for every $f\in\mathrm{Hom}\,(A,B)$.

(2) $\pi(H)$ is a characteristic subgroup in $A$.
\end{lemma}

\Pf (1) It is obvious that the endomorphism
$\left(\begin{array}{cc} 1 & 0 \\ f & 1\end{array}\right)$
is an automorphism of $G$, hence, for any $x=a+b\in H$ with $a\in A$, $b\in B$, we have $f(x)=x+f(a)$. Thus, $f(a)=f\pi(x)\in H$, as required.

(2) If $\alpha\in\mathrm{Aut}(A)$, then one checks that
$\left(\begin{array}{cc} \alpha & 0 \\ 0 & 1\end{array}\right)\in\mathrm{Aut}(G)$, hence,
for any $a+b\in H$ with $a\in A$, $b\in B$, we have $\alpha(a)+b\in H$. So, $a-\alpha(a)\in H$ and, therefore,
$\alpha(a)=a-(a-\alpha(a))\in\pi(H)$, as required.
\fine

\begin{proposition}\label{02} (1) Let $G=\bigoplus_{i\in I} G_i$, where, for each direct summand $G_i$, there is an other summand $G_j$ with $j\in I\setminus \{i\}$ such that $G_i\cong G_j$. Then every characteristic subgroup of
$G$ is fully invariant.

\medskip

(2) Let $A=B\oplus G$, where $G$ is such a decomposable group that, for every its direct summand $C$, it decomposes as $G=C\oplus K$ and in the complement $K$ there exists a direct summand isomorphic to $C$. Then, if $\pi:A\to G$ is projection and $H$ is a characteristic subgroup in $A$, it follows that $\pi(H)=H\cap G$ and $H\cap G$ is fully invariant in $G$.
\end{proposition}

\Pf (1) In view of the isomorphism between $G_i$ and $G_j$ and consulting with Lemma~\ref{01}, we deduce that not only that $\pi_i(H)\leq H$ for the projection $\pi_i: G\to G_i$, but also that $\varphi\pi_i(H)\leq H$ for every
endomorphism $\varphi$ of $G_i$ and, besides, that $f\pi_i(H)\leq H$ for every $f\in\mathrm{Hom}\,\bigl(G_i,\bigoplus_{j\neq i}G_j\bigr)$. This surely means that the subgroup $H$ is fully invariant, as claimed.

(2) follows directly from (1).
\fine

We can now say slightly more in the case of arbitrary direct sums.

\begin{lemma}\label{matrices} Let $G=\bigoplus_{i\in I}G_i$ and $\pi_i: G\to G_i$ be the corresponding projections, where $I$ is an arbitrary index set. If $H$ is a characteristic subgroup in $G$, then both the subgroups $\underline{H}=\bigoplus_{i\in I}(H\cap G_i)$ and $\overline{H}=\bigoplus_{i\in I}\pi_i(H)$ also are characteristic subgroups in $G$.
\end{lemma}

\Pf Let $x\in H\cap G_i$ and $\alpha\in\mathrm{Aut}(G)$. Since $\alpha(x)$ can be embedded in a sum of finite numbers
of direct components $G_i$, we may assume that the set $I$ is finite, thus writing $G=G_1\oplus\dots\oplus G_n$ for some $n\geq 2$ and $x\in G_1$. Then, mimicking \cite{F},\cite{F1}, we may identify $\alpha$ as the $n\times n$ matrix

$$\alpha=\left(\begin{array}{cccc} \alpha_{11} & \alpha_{12} & \dots & \alpha_{1n}\\
\alpha_{21} & \alpha_{22} & \dots & \alpha_{2n}\\
\dots & \dots & \dots & \dots\\
\alpha_{n1} & \alpha_{n2} & \dots & \alpha_{nn}\end{array}\right)$$
and so we write
$\alpha(x)=\alpha_{11}(x)+(\alpha_{21}(x)+\dots+\alpha_{n1}(x))\in H$.
We have $\alpha_{11}(x)\in G_1$ and hence Lemma~\ref{01} (1) yields $\alpha_{i1}(x)\in H\cap G_i$, where $i=2,\dots,n$, so that $\alpha_{11}(x)\in H\cap G_1$. Consequently, $\underline{H}$ is characteristic in $G$, as claimed.

If now $x_1\in\pi_1(H)$, then $x_1+\dots+x_n=h$ for some $h\in H$ and $x_i\in G_i$ ($i=2,\dots, n$).
Therefore,
$$\alpha(h)=\alpha_{11}(x_1)+\alpha_{22}(x_2)+\dots+\alpha_{nn}(x_n)+
\sum_{i\neq j}\alpha_{ij}(x_{j}),$$
where $\alpha(h)\in H$ and again applying Lemma~\ref{01} (1) it must be that $\alpha_{ij}(x_{j})\in H\cap G_i$
whenever $i\neq j$.
So, we get
$$\pi_1\alpha(h)=\alpha_{11}(x_1)+\alpha_{12}(x_2)+\dots+\alpha_{1n}(x_n)\in \pi_1(H).$$
Since $\alpha_{12}(x_2)+\dots+\alpha_{1n}(x_n)\in \pi_1(H)$, one has that $\alpha_{11}(x_1)=\pi_1\alpha(x_1)\in\pi_1(H)$ and $\pi_i\alpha(x_1)\in H\cap G_i\leq\pi_i(H)$, provided
$i=2,\dots,n$, i.e., $\alpha(x_1)\in\pi_1(H)\oplus\dots\oplus\pi_n(H)$. Thus, $\pi_1(H)\oplus\dots\oplus\pi_n(H)$
is characteristic in $G$, as asserted.
\fine

\begin{remark} Note that if $G=A\oplus B$ and $\{0\}\neq H\leq G$ is a characteristic subgroup with $H\nleq B$, then
$H\cap A\neq \{0\}$ presuming the condition $A[2]=\{0\}$. Indeed, there exists $h\in H$ with $h=a+b$, where $a\neq 0$. Since $h'=a-b$ lies also in $H$, we have $0\neq h+h'=2a\in H\cap A$.

However, it is constructed in Proposition~\ref{10}, listed below, a $2$-group $G$ having a direct summand $\langle a\rangle\neq \{0\}$ and a characteristic subgroup $\{0\}\neq H\leq G$ with $H\cap \langle a\rangle=\{0\}$, as needed.
\end{remark}

Notice that, in view of Lemma~\ref{matrices}, if $G$ is a direct sum of cyclic groups, writing $G=\bigoplus_{n\geq 1}G_n$, where $G_n$ is a direct sum of cyclic groups of order $p^n$ and $\pi_n:G\to G_n$ are the corresponding projections, and $\overline{H}=\bigoplus_{n\geq 1}\pi_n(H)$ for some characteristic subgroup $H$ of $G$, then $\overline{H}$ is a fully invariant subgroup such that $\overline{H}\leq F$ for every fully invariant subgroup $F$ with $H\leq F$. In fact, exploiting Proposition~\ref{10} quoted below, each subgroup $\pi_n(H)$ is fully invariant in $G_n$ and thus it follows from Lemma~\ref{01} that $\overline{H}$ is really a fully invariant subgroup
in $G$. Furthermore, since $\pi_n(H)\leq F$ for every fully invariant subgroup $F$ with $H\leq F$, it can be deduced that $F$ contains $\overline{H}$, as expected.

Note also that each fully invariant subgroup $H$ of $G$ is closed in $G$, that is, the quotient-group $G/H$ is separable. Indeed, this fact follows from the observation that $H=\bigoplus_{n\geq 1}(H\cap G_n)$, where the factor-groups $G_n/(H\cap G_n)$ are also a direct sum of cyclic groups. However, if we assume that $H$ is only characteristic in $G$, this property may fail as the next example manifestly shows (in it $G$ must be a $2$-group).

\begin{example}
Let $G=\langle a_1\rangle\oplus\langle a_2\rangle\oplus\dots$, where $o(a_i)=2^{2i-1}$, $h_i=a_1+2^ia_{i+1}$, $H=\langle h_i\,|\,i\geq 1\rangle$. Then $H$ is not a closed characteristic subgroup in $G$.
\end{example}

\Pf Write $a_1=h_i-2^ia_{i+1}$ for every $i\geq 1$, i.e., $\overline{0}\neq a_1+H\in (G/H)^1$, which means that $H$ is not closed in $G$. Moreover, the subgroup $\langle h_i\rangle$ is characteristic in $\langle a_1\rangle\oplus\langle a_{i+1}\rangle$ as it is generated by all its elements of the kind $o(x)=2^{i+1}$, $h(x)=0$ and $h(2x)=i+1$. And since one easily checks that the two inequalities
$$f(\langle 2^{i+1}a_{i+2}\rangle)\leq 2^{i+1}a_{i+1}\leq\langle h_{i}\rangle,
\varphi(2^{i}a_{i+1})\leq\langle 2^{i+2}a_{i+2}\rangle\leq\langle h_{i+1}\rangle$$
hold for each element $f\in\mathrm{Hom}\bigl(\bigoplus_{j\geq i+2}\langle a_j\rangle,\bigoplus_{1\leq j\leq i+1}\langle a_j\rangle\bigr)$ and each element $\varphi\in\mathrm{Hom}\bigl(\bigoplus_{1\leq j\leq i+1}\langle a_j\rangle,\bigoplus_{j\geq i+2}\langle a_j\rangle\bigr)$, we then can infer that $H$ is characteristic in $G$, as promised.
\fine

The next series of three propositions are quite useful in the sequel.

\begin{proposition}\label{03} Every bounded $p$-group, divisible $p$-group and divisible torsion-free group is not a weakly IC-group.
\end{proposition}

\Pf Let $G$ be a bounded $p$-group, say $p^{k+1}G=0$ with $p^{k}G\neq 0$, and assume the contrary that $H\cong G$ is characteristic in $G$. If $o(g)=p^k$ for some $g\in H$, then we have $G=\langle g\rangle\oplus A$
for some $A\leq G$. Since $\mathrm{Hom}\,(\langle g\rangle,A)\langle g\rangle=A$, with Lemma~\ref{01} at hand we derive that $H=G$, a contradiction. Furthermore, as any divisible subgroup is a direct summand of a given group, the two remained cases can be proved similarly.
\fine

\begin{remark} It follows from Proposition~\ref{03} that all bounded groups and divisible groups are not weakly IC-groups. Moreover, it is easy to see that a non-reduced group is a weakly IC-group if, and only if, so is its reduced part. And so, we hereafter will consider only {\bf reduced} groups.

On the same vein, since each torsion-free weakly IC-group $G$ must be non-divisible, one sees that $nG\cong G$ but $nG\neq G$ for some natural number, so for the torsion-free weakly IC-group $G$ we shall assume that it contains own characteristic subgroups different from $nG$ which are isomorphic to $G$ (notice that this observation held too in \cite{GN3} for the case of IF-groups).
\end{remark}

\begin{proposition}\label{04} An arbitrary non-reduced torsion-free group is a weakly IC-group if, and only if, its reduced part is a weakly IC-group.
\end{proposition}

\Pf "{\bf Necessity.}" Write $G = R\oplus D$, where $R$ is the reduced part and $D$ is the divisible part of the group $G$, and suppose $\pi: G\to R$ is the corresponding projection. Letting $H$ be a characteristic subgroup in
$G$, if $\pi(H) = \{0\}$, then $H\leq D$ which implies that $H = D$. So, if $H\cong G$, then $\pi(H)\neq \{0\}$. However, it follows from Lemma~\ref{01} that $D\leq H$. Thus, $H = (H\cap R)\oplus D$. Therefore, if $H\cong G$, then $H\cap R\cong R$, as required.

The "{\bf Sufficiency}" is obvious.
\fine

Note that, according to \cite[Corrolary 6]{GN3}, any fully invariant subgroup of a homogeneous fully transitive torsion-free group $G$, isomorphic to the whole group $G$, is of the form of $nG$ for some natural number $n$. As homogeneous transitive torsion-free groups are always fully transitive, we obtain:

\begin{corollary}\label{06} Homogeneous fully transitive torsion-free groups are not weakly IC-groups.
\end{corollary}

We, therefore, obtain the following construction: chosen a homogeneous separable torsion-free group of an arbitrary infinite power $\mathfrak{m}$, we thus can get the existence of non-weakly IC-groups.

And since in each torsion-free group $G$ whose ring of endomorphisms is isomorphic to $\mathbb{Z}$, any endomorphic image is a subgroup $nG$ for some natural number $n$, we extract:

\begin{proposition}\label{07} The torsion-free group whose endomorphism ring is isomorphic $\mathbb{Z}$ is not a weakly IC-group.
\end{proposition}

Note that, for any natural number $n$ and any infinite cardinal number $\mathfrak{m}$, which is the smaller first powerfully unattainable number, there exist such torsion-free groups of rank $n$ and torsion-free groups of power $\mathfrak{m}$ (see \cite[\S~ 88, Excercise 8; \S~89, Excercise 1]{F}).

\begin{example}\label{08} There exist indecomposable torsion-free groups which are a weakly IC-group.
\end{example}

\Pf We shall exhibit two constructions as follows:

\medskip

(1) Take the group $G$ with $\mathrm{E}(G) = \mathbb{Z}+\sqrt{-5}\mathbb{Z}$. Then, in this group, all endomorphisms are injective and all subgroups are characteristic. So, $G\cong\sqrt{-5}G$. Furthermore, assume that $\sqrt{-5}G = kG$ for some natural number $k$. Consequently, $5G = k^2G$. Since the numbers $5$ and $k$ are mutual simple, one checks that $kG = G$ and hence $k = 1$. However, the inequality $\sqrt{-5}G\neq G$ is true, that is the desired contradiction. Note that such a group $G$ is possible to be chosen as given (e.g., of rank $4$, it is exhibited in \cite[\S~110, Excercise 8]{F}).

(2) Let $\hat{\mathbb{Z}}=\prod_{p\in\mathbb{P}}\hat{\mathbb{Z}}_p$, where $\mathbb{P}$ is the set of all prime numbers, and $\hat{\mathbb{Z}}_p$ is the ring of $p$-adic integers. In addition, $\pi_p: \hat{\mathbb{Z}}\to \hat{\mathbb{Z}}_p$ are the corresponding projections, $1=(\dots,1_p,\dots)$ is the unit of $\hat{\mathbb{Z}}$,
and $\xi\in \hat{\mathbb{Z}}$ is such that $\pi_p(\xi)=p\xi_p\in\hat{\mathbb{Z}}_p$, where $\xi_p$ is some transcendental element on $\mathbb{Z}$. Further, we define $S=\mathbb{Z}[\xi]$ to be the subring in $\hat{\mathbb{Z}}$ induced by $\{1,\xi\}$, and $R=S^+_*$ is a pure hull of $S^+$ in $\hat{\mathbb{Z}}$. Thus $R$ will also be a subring in $\hat{\mathbb{Z}}$. It is readily checked that each element of $R$ is an integer multiple of some $a_0\xi^n+a_1\xi^{n-1}+\dots+a_n$ with $a_i\in\mathbb{Z}$, so the only invertible elements in $R$ are $\pm 1$. Since
$G=R^+$ is an $E$-group as being a pure subgroup of the $E$-group $\hat{\mathbb{Z}}$, one derives that any subgroup in $G$ is really characteristic. But we have $G\cong\xi^2 G$, and because the type of each element from $\xi^2 G$ is surely strictly more than the type of $\xi$, we conclude that $\xi^2 G\neq kG$ for all $k\in\mathbb{N}$, as required.

Note that the rings for which each element is an integer multiple of invertible are identified as {\it strongly homogeneous rings}.
\fine

If $G$ is an arbitrary weakly IC-group, then a direct summand of $G$ cannot be a weakly IC-group. Indeed, such a possibility occurs by taking a bounded direct summand of $G$ (see, for instance, Proposition~\ref{03}).

\medskip

We shall now deal with primary groups by starting with the following example which manifestly illustrates that the class of weakly IC-groups is {\it not} closed under the formation of direct summands.

\begin{example}\label{direct2} There exists such a $p$-primary weakly IC-group $G$ that some unbounded direct summand of $G$ is not a weakly IC-group.
\end{example}

\Pf By virtue of \cite[Theorem~19 and Corrolary~21]{GN2}, if $G$ is a direct sum of cyclic $p$-groups such that $f_n(G)=\gamma$, $n\in\mathbb{N}_0$, for some cardinal $\gamma$, then $G$ is a weakly IF-group. Letting $\gamma\geq\aleph_0$, then Proposition~\ref{02} (1) tells us that every characteristic subgroup of $G$ is fully invariant, so that $G$ is a weakly IC-group. Suppose $A$ is such a direct summand of $G$ that its Ulm-Kaplansky invariants have the properties $f_0(A)\geq 2$, provided $f_0(A)\neq 0$, and all $f_n(A)$ form an increasing sequence of non-negative integers. Thus we can apply \cite[Corrolary~18]{GN2} to get that $A$ is not a weakly $IF$-group, and hence the posed restrictions on $f_n(A)$ guarantee that $A$ is also not an weakly IC-group.
\fine

Observe that the following technicality is similar but independent to Theorem 2.8 from \cite{B}.

\begin{lemma}\label{09} Let $B = \bigoplus_{k\geq 1} B_k$, where $B_k = \bigoplus\mathbb{Z}(p^k)$ and $\pi_k: B\to B_k$ are the corresponding projections, $H$ is a characteristic subgroup of $B$. Then $\pi_k(H) = p^{n_k}B_k$, where:

(1) $n_k\leq k$ for all $k\geq 1$;

(2) $n_k\leq n_{k+r}\leq n_{k}+r$ for all $k,r\geq 1$;

(3) if $n_k<k$, then $p^{n_k+r}B_{k+r}\leq H$ for all $k,r\geq 1$;

(4) if $r(B_k)\geq 2$, then $\pi_k(H)\leq H$.
\end{lemma}

\Pf Assume that $\pi_k(H)\neq \{0\}$. If $r(B_k) = 1$, then $\pi_k(H) = p^{n_k}B_k$ for some $n_k < k$. Let $r(B_k) > 1$, $B_k = \bigoplus_{i\in I}\langle a_i\rangle$ and $\theta_i: B\to\langle a_i\rangle$ are the corresponding projections. Thus, $\theta_i(H)=\langle p^{n_{k_i}}a_i\rangle$. Let $n_k = \min\{n_{k_i}\,|\,i\in I\}$. If $f: a_k \to  a_i$ is an isomorphism, then $p^{n_k}a_i = f(p^{n_k}a_k)\in H$ in virtue of Lemma~\ref{01}. So, in this case, $\pi_k(H) = p^{n_k}B_k\leq H$, i.e., points (1) and (4) are fulfilled.

(2) Let $x\in p^{n_k}B_k\setminus p^{n_k+1}B_k$. Thus $o(x) = p^{k-n_k}$, so $f(x)\in p^{n_{k}+r}B_{k+r}\leq H$ for every $f\in\mathrm{Hom}\,(B_k,B_{k+r})$. Therefore, $n_{k+r}\leq n_{k}+r$. If $y\in H\cap p^{n_{k+r}}B_{k+r}$, where $n_{k+r}<n_k$, then $\varphi(y)\in p^{n_{k+r}}B_{k}$ for some $\varphi\in\mathrm{Hom}\,(B_{k+r},B_k)$, because $\varphi(y)\in H$ contradicts $\pi_k(H) = p^{n_k}B_k$. Finally, the inequalities $n_k\leq n_{k+r}\leq n_{k}+r$ hold.

(3) The proof is analogous to the previous point. In fact, if $B_k = \bigoplus_{i\in I}\langle a_i\rangle$, then for every $x\in p^{n_k+r}B_{k+r}$ such that $n_k<k$ there is $f\in\mathrm{Hom}\,(B_{n_k},B_{k+r})$
with $f(p^{n_k}a_i)=x$. Applying now Lemma~\ref{01}, we conclude that $x\in H$, i.e., $p^{n_k+r}B_{k+r}\leq H$, as wanted.
\fine

The next statement is often attributed to Kaplansky who showed it for modules (cf. \cite[Theorem~26,~p.~63]{Kap}).

\begin{proposition}\label{10} If $G = \bigoplus_{i\geq 1}G_i$ is a $2$-group, where $G_i=\bigoplus \mathbb{Z}(2^i)$ is a direct sum of cyclic groups, then each characteristic subgroup of $G$ is fully invariant if, and only if, $G$ has at most two Ulm invariants equal to one, and if it has exactly two, they correspond to successive ordinals.
\end{proposition}

\Pf "$\Rightarrow$". Assume that $f_{k-1}(G) = f_{k+r-1}(G) = 1$, where $r > 1$, and $\langle a\rangle\oplus\langle b\rangle$ is the corresponding direct summand of $G$ with $o(a) = 2^k$, $o(b) = 2^{k+r}$. Suppose $h = 2^{k-1}a+2^{k+r-2}b$, $H=\langle h\rangle = \{0, 2^{k-1}a\pm 2^{k+r-2}b, 2^{k+r-1}b\}$ and
$$\underline{H}=H\oplus\left(\bigoplus_{k+1\leq i\leq k+r-1}G_i\right)[2]\oplus \left(\bigoplus_{i\geq k+r+1}G_i\right)[2^2].$$
Thus, one verifies that the subgroup $\underline{H}$ has to be characteristic in $G$. Indeed, $H$ is characteristic in $\langle a\rangle\oplus\langle b\rangle$, $f(\underline{H}) = \{0\}$ for $f\in\mathrm{Hom}\,(\bigoplus_{i\geq k}G_i,\bigoplus_{i<k}G_i)$, $\varphi(H)\leq\underline{H}$ for $\varphi\in\mathrm{Hom}\,\bigl(\langle a\rangle\oplus\langle b\rangle,\bigoplus_{i>k,\, i\neq k+r}G_i\bigr)$ and $\bigl(\bigoplus_{k+1\leq i\leq k+r-1}G_i\bigr)[2]\oplus\bigl(\bigoplus_{i\geq k+r+1}G_i\bigr)[2^2]$ is fully invariant in $\bigl(\bigoplus_{k+1\leq i\leq k+r-1}G_i\bigr)\oplus\bigl(\bigoplus_{i\geq k+r+1}G_i\bigr)$. But $H$ is not fully invariant since for the projection $\alpha$ of $G$ on $\langle a\rangle$, we have $\alpha(H)=\langle 2^{k-1}a\rangle\nleq H$, as expected.

"$\Leftarrow$". Suppose that $G_k=\langle a\rangle$, $G_{k+1} = \langle b\rangle$ and $r(G_i)\geq 2$ with $i\neq k, k+1$. Let $\pi: G\to G_k$, $\theta: G\to G_{k+1}$ be the corresponding projections and let $H\leq G$ be a  characteristic subgroup of $G$. It follows with the aid of Lemma~\ref{01} that if $h = x_i+y$, where $x_i\in G_i$ for some $i\neq k, k+1$, then $x_i\in H$ and $\beta(H)\leq H$, where $\beta$ is the projection of $G$ on
$\bigoplus_{j\neq k,\,k+1}G_j$, and $\beta(H)$ is fully invariant in $G$. Letting $\pi(H) = p^{n_{k}}a$, $\theta(H) = p^{n_{k+1}}b$, then $p^{n_{k}+1}b\in H$. If $p^{n_{k}+1}= p^{n_{k+1}}$, then $p^{n_{k}}a\in H$
since in this case $p^{n_{k}}a+y\in H$ for some $y\in H$, as it was noted above, and so $p^{n_{k}}a\in H$. One inspects now that $H$ is really fully invariant in $G$. If, however, $p^{n_{k}+1}> p^{n_{k+1}}$, then
$p^{n_{k+1}} \leq p^{n_{k}}$, but in this situation the relation $p^{n_{k+1}}a\in H$ is possible only if $p^{n_{k+1}} = p^{n_{k}}$, hence $\pi(H)\leq H$ which yields $\theta(H)\leq H$, i.e., $H$ is fully invariant. In the remaining case when $p^{n_{k}+1}< p^{n_{k+1}}$, it can be inferred that $p^{n_{k}+1}b\in H$ that is impossible.
\fine

We can, alternatively, extend the last result to the following one.

\begin{proposition}\label{101} Given $G$ is a separable $2$-group, then each characteristic subgroup of $G$ is fully invariant if, and only if, $G$ has at most two Ulm invariants equal to one, and if it has exactly two, they correspond to successive ordinals, i.e., if, and only if, the corresponding properties are satisfied by its basic subgroup.
\end{proposition}

\Pf "$\Rightarrow$". In the same manner, as showed above in the preceding proposition, if $B=\bigoplus_{i\geq 1}B_i$, where $B_i=\bigoplus \mathbb{Z}(2^i)$, is a basic subgroup in $G$, then we may write
$G=\left(\bigoplus_{i=1}^{n}B_i\right)\oplus\left(B_n^*+2^nG\right)$ for any $n\in\mathbb{N}$, where
$B_n^*=\bigoplus_{i\geq n+1}B_i$.

Let us now $f_{k-1}(G) = f_{k+r-1}(G) = 1$, where $r > 1$, and $\langle a\rangle\oplus\langle b\rangle$ is the corresponding direct summand of $G$ with $o(a) = 2^k$, $o(b) = 2^{k+r}$. Suppose $h = 2^{k-1}a+2^{k+r-2}b$, $H=\langle h\rangle = \{0, 2^{k-1}a\pm 2^{k+r-2}b, 2^{k+r-1}b\}$ and
$$\underline{H}=H\oplus\left(\bigoplus_{k+1\leq i\leq k+r-1}B_i\right)[2]\oplus \bigl(B_{k+r}^*+2^{k+r}G\bigr)[2^2].$$
Thus, one inspects that $\underline{H}$ is a characteristic, but not fully invariant subgroup in $G$. Note that the separability in this case is not used.

"$\Leftarrow$". Let $H$ be a characteristic subgroup in $G$, and let $x\in H$ and $y=f(x)$ for some $f\in\mathrm{E}(G)$. Since $G$ is separable, we know that the elements $x$ and $f(x)$ can be embedded in some basic subgroup, say $B$, of $G$. Consulting with Lemma~\ref{13} stated and proved below, one finds that $H\cap B$ is characteristic in $B$, so it is fully invariant. And since $\mathrm{h}^B(x)\leq\mathrm{h}^B(f(x))$ as well as $B$ is fully transitive, it must be that $\varphi(x)=f(x)$ for some $\varphi\in\mathrm{E}(B)$, whence $f(x)\in H$, as required.
\fine

The next technical claim is rather curious and very helpful for our objectives.

\begin{lemma} If $G$ is a $p$-group and $H$ is its characteristic subgroup such that the quotient-group $H/(H\cap p^{\omega}G)$ is unbounded, then $p^{\omega}G\leq H$.
\end{lemma}

\Pf If $B=\bigoplus_{i\geq 1}B_i$, where $B_i=\bigoplus \mathbb{Z}(p^i)$, is a basic subgroup in $G$, then we can write $G=\left(\bigoplus_{i=1}^{n}B_i\right)\oplus\left(B_n^*+p^nG\right)$ for each $n\in\mathbb{N}$, where
$B_n^*=\bigoplus_{i\geq n+1}B_i$. Let $\pi_n: G\to B_1\oplus\dots\oplus B_n$ be the corresponding projections.
Since the factor-group $H/(H\cap p^{\omega}G)$ is unbounded, for every $s\in\mathbb{N}$, there exist elements $x=p^sh$, $h\in H$ with $x\notin p^{\omega}G$. Note that $p^s\pi_{n_s}(h)\neq 0$ for some $n_s>s$ as for otherwise we will have $x\in B_{n_s}^*+p^{n_s}G$ for each $n_s>s$. If now $o(h)=p^{m_s}$, then $m_s\geq s$. By assumption, we get $x\in B_{k+m_s}^*+p^{k+m_s}G$ for all $k\geq 1$, so $x\in p^{k+1}G$, i.e., $x\in p^{\omega}G$, thus contradicting the
choice of the element $x$. Therefore, for every natural $t$, there will exist $n_t\in\mathbb{N}$ and $x\in H$ such that $o(\pi_{n_t}(x))\geq t$. But since $U(y)\geq U(\pi_{n_t}(x))$ for each $y\in p^{\omega}G$ with $o(y)\leq t$, we derive that $f(\pi_{n_t}(x))$ for some $f\in\mathrm{Hom}(B_1\oplus\dots\oplus B_{n_t},B_{n_t}^*+p^{n_t}G)$. Thus, it follows from Lemma~\ref{01} that $y\in H$. Finally, $p^{\omega}G\leq H$, as formulated.
\fine

The following assertion is closely related to the previous one.

\begin{proposition} If $G$ is a $p$-group such that either $f_n(G)=0$ or $f_n(G)\geq 2$ for all $n\in\mathbb{N}\cup \{0\}$, then, for every characteristic subgroup $H$ of $G$, the subgroup $H'=H+p^{\omega}G$ is fully invariant in $G$. In particular, if $H/(H\cap p^{\omega}G)$ is unbounded, then $H$ is fully invariant in $G$.
\end{proposition}

\Pf Having in mind Proposition~\ref{101}, the quotient $H'/p^{\omega}G$ is fully invariant in $G/p^{\omega}G$, and consequently $H'$ is fully invariant in $G$, as stated. The second part is now obvious by using the argumentation from the previous statement.
\fine

To avoid any further misunderstanding and confusion, let us recall that $p^nG[p^m]$ just means $(p^nG)[p^m]$ written without brackets for simpleness.

\medskip

The next four technical assertions are pivotal for obtaining our chief results. They somewhat expand closely related results from \cite{F} and \cite{F1} (see also the existing bibliography therewith).

\begin{lemma}\label{11} If $B$ is a basic subgroup of a reduced $p$-group $G$, then the following two points are valid:

(1) $G[p^m]\leq B+ p^nG[p^m]$ for all $m,n\in \mathbb{N}$;

(2) $G = B+\sum_{i\geq 1}p^{n_i}G[p^{m_i}]$ for any sequences of natural numbers $1\leq m_1< m_2<\dots$ and for all set of natural numbers $\{n_i\}$.
\end{lemma}

\Pf (1) Write $G=B+p^nG$ for all $n\in \mathbb{N}$. If now $p^mg=0$ for some $g\in G$, then for $g=b+p^ng'$, $g'\in G$, we have $p^mb+p^{n+m}g'=0$. Since $B$ is pure in $G$, it follows that $p^mb=p^{n+m}b'$ for some $b'\in B$. So, $p^{n+m}g'=-p^{n+m}b'$ and, consequently, $g'+b'\in G[p^{n+m}]$. Thus, $p^n(g'+b')=g''$ for some $g''\in p^nG[p^m]$, whence $g=b-p^nb'+g''\in B+p^nG[p^m]$, as desired.

(2) follows from (1) and the fact that $G=\sum_{i\geq 1}G[p^{m_i}]$ for every sequences of natural numbers $1\leq m_1< m_2<\cdots$.
\fine

\begin{lemma}\label{12} If $H$ is a non-bounded characteristic subgroup of a separable $p$-group $G$, then $G=H+B$ for every basic subgroup $B$ of $G$.
\end{lemma}

\Pf Employing Lemma~\ref{11} it is enough to show that for the sequences of natural numbers $1\leq m_1<m_2<\cdots$ the inclusion $\sum_{i\geq 1}p^{n_i}G[p^{m_i}]\leq H$ holds for some set of natural numbers $\{n_i\}$. To show that, for every $n\geq 1$ we have $G=(B_1\oplus\dots\oplus B_{n+1})\oplus(B_{n+1}^*+p^{n+1}G)$, where $B_{n+1}^*=\bigoplus_{i\geq n+2}B_i$. Let $\pi_n:G\to B_n$ are the corresponding projections. Let $h_1\in H$ and $p^{n+1}h_1=0$, $p^{n}h_1\neq 0$ for some $n\in\mathbb{N}$. Note that $\pi_{n+k}(p^nh_1)\neq 0$ for some $k\geq 1$. In fact, if $\pi_{n+k}(p^nh_1)= 0$ for every $k\geq 1$, then $p^nh_1\in (B_{n+k}^*+p^{n+k}G)[p]$.
However, since $p^{n+k}G$ is essential in $B_{n+k}^*+p^{n+k}G$, we infer that $(B_{n+k}^*+p^{n+k}G)[p]\leq p^{n+k}G$. Therefore, $0\neq p^nh_1\in\bigcap_{k\geq 1}p^{n+k}G=\{0\}$ -- contradiction. So, $0\neq \pi_{n+k}(p^{n_1}h)= p^{n}\pi_{n+k}(h_1)\in B_{n+k}$ for some $k\geq 1$. Suppose now $p^m_1=o(\pi_{n+k}(h_1))$. Since $U(\pi_{n+k}(h_1))\leq U(z)$ for each $z\in p^{n+k-m_1}(B_{n+k}^*+p^{n+k}G)[p^m_1]$, Lemma\ref{01} assures that $p^{n+k-m_1}(B_{n+k}^*+p^{n+k}G)[p^m_1]\leq H$, and since $p^{n+k-m_1}G= p^{n+k-m_1}(B_{n+k-m_1}^*+p^{n+k-m_1}G)$, it must be that $p^{n_1} G[p^{m_1}]\leq H$, where $n_1=n+k-m_1$.

Furthermore, since the subgroup $H$ is unbounded, there will exist a natural number $m_2>m_1$ and some positive integer $n_2$ satisfying the relation $p^{n_2} G[p^{m_2}]\leq H$, as wanted.
\fine

\begin{lemma}\label{13} If $H$ is an unbounded characteristic subgroup of a separable $p$-group $G$ and $B\leq G$ is a basic, then $H\cap B$ is basic in $H$.
\end{lemma}

\Pf Our first goal is to demonstrate that $H\cap B$ is pure in $H$. Indeed, let $px=b$ for some $x\in H$, $b\in H\cap B$. Since $B$ is pure in $G$, we may write $b=py$ for some $y\in B$. If $h(x)=h(y)$, then $U(x)=U(y)$,
and since $G$ is transitive we write $\alpha(x)=y$ for some $\alpha\in\mathrm{Aut}(G)$; so, $y\in H\cap B$. Assuming that $m=h(x)<h(y)$, we then have $h(px)>m+1$ and thus $f_m(G)\neq 0$; consequently, $B_{m+1}\neq \{0\}$.
Letting $z\in B_{m+1}[p]$, since $h(z)=m$, we deduce $h(y+z)=h(z)=h(x)$, and so $U(x)=U(y+z)$; therefore, as above, $y+z\in H\cap B$, where $p(y+z)=b$. Let us now $h(x)>h(y)=m$. So, $B_{m+1}\neq \{0\}$ and there exists $z\in B_{m+1}[p]$ such that $h(y-z)>h(y)$. If, however, $h(x)>h(y-z)$, we similarly can choose $z_1\in B[p]$ with $h(y-z-z_1)>h(y-z)$. Furthermore, after a finite number of steps, we get either $h(x)=h(y+u)$, or $h(x)<h(y+u)$
for some $u\in B[p]$. In both cases, as shown above, the equation $px=b$ is resolvable in $H\cap B$. But since $B+pH$ also equals to $G$, it follows that $H\cap p^nH=p^n(H\cap B)$ for every $n$, i.e., $H\cap B$ is pure in $H$. And since $(B+H)/B\cong H/(H\cap B)$, where $B+H=G$, the subgroup $H\cap B$ is dense in $H$, i.e., $H\cap B$ is basic in $H$, as claimed.
\fine

\begin{lemma}\label{14} If $H$ is a characteristic subgroup of a separable $p$-group $G$ and $B$ is basic in $G$,  then $H\cap B$ is characteristic in $B$.
\end{lemma}

\Pf For every $n$, we have that $G=(B_1\oplus\dots\oplus B_n)\oplus(B_n^*+p^nG)$. Let $\pi_n:G\to B_n$ are the corresponding projections. If $H\cap B\neq \{0\}$, then $\pi_n(H)\neq \{0\}$ for some $n$. If $r(B_n)>1$, then consulting with Proposition~\ref{02} we inspect that $\pi_n(H)\leq H$ is fully invariant in $B_n$. Moreover, Lemma~\ref{01} tells us that $f(H\cap B_n),\varphi(H\cap B_n)\leq H$ for each $f\in\mathrm{Hom}\,(B_n, B_1\oplus\dots\oplus B_{n-1})$ and $\varphi\in\mathrm{Hom}\,(B_n, B_n^*+p^nG)$.

Let us now $x\in H\cap B$ and $x=x_1+\dots +x_m$, where $x_j=\pi_{i_j}(x)\in B_{i_j}$. In case that $r(B_{i_j})>1$, we can observe that $\pi_{i_j}(H)\leq H$ and that $\pi_{i_j}(H)$ is fully invariant in $B_{i_j}$, so we may  assume that $r(B_{i_j})=1$ for all $j=1,\dots,m$. Given $\alpha\in\mathrm{Aut}(B)$, we note that $\pi_{i_j}\alpha$ is an automorphism for $B_{i_j}$ for all $j=1,\dots,m$. In fact, if $\pi_{i_j}\alpha(y)=0$ for some $0\neq y\in B_{i_j}[p]$, then $\pi_k\alpha(y)=0$ for every $k<i_j$ and so $0\neq\alpha y\in (B_{i_j}^*)[p]$. But, in this case, $U(\alpha y)>U(y)$ which is obviously imposable for getting automorphisms. So, $\ker(\pi_{i_j}\alpha)=0$, and in view of the finite number of the $B_{i_j}$s, it follows that $(\pi_{i_j}\alpha) \!\upharpoonright\! B_{i_j}$ is an automorphism for all $j=1,\dots,m$. Consequently,
$$\mu=(\pi_{i_1}+\dots+\pi_{i_m})\alpha \!\upharpoonright\! \bigoplus_{j=1}^mB_{i_j}$$
is an automorphism of $\bigoplus_{j=1}^mB_{i_j}$. Since $\bigoplus_{j=1}^mB_{i_j}$ is a direct summand of $B$, $\mu$ can be considered as an extension of the automorphism of $B$, so $(\pi_{i_1}+\dots+\pi_{i_m})\alpha(x)\in H\cap B$, as promised, because it follows from Lemma~\ref{01} that $(1_B-(\pi_{i_1}+\dots+\pi_{i_m}))\alpha(x)\in H\cap B$, as asked for.
\fine

We are now in a position to prove one of our main results. It is worthwhile noticing that the following theorem is similar to Theorem 14 in \cite{GN1} proved for IF-groups.

\begin{theorem}\label{15} The separable $p$-group $G$ is not a weakly IC-group if some of its basic subgroups $B$ is not a weakly IC-group.
\end{theorem}

\Pf If $G$ is bounded, then its basic subgroups coincide with $G$, hence Proposition~\ref{03} applies to get that $G$ is not a weakly IC-group, as asserted.

Let us now $G$ be unbounded, and assume in a way of contradiction that $G$ is a weakly IC-group. Then there exists a characteristic subgroup $H\lneqq G$ such that $H\cong G$. Since $H$ is unbounded, one derives by Lemma~\ref{13} that $H\cap B$ is basic subgroup in $H$. If $H\cap B=B$, then $B\leq H$, and since $H+B=G$, the application of Lemma~\ref{12} guarantees that $H=G$, that is a contrary to our initial choice. So, $H\cap B\neq B$.
Since $H\cong G$, its basic subgroups are also isomorphic, say $H\cap B\cong B$ (see, for instance, \cite{F}). However, Lemma~\ref{14} ensures that $H\cap B$ is characteristic in $B$, which is an absurd, so we get the pursued contradiction.
\fine

The next technical claim is the key in proving the opposite direction of the last theorem.

\begin{lemma}\label{16} (1) If $H$ is an unbounded characteristic subgroup in a fixed basic subgroup $B$ of a torsion-complete $p$-group $G$, then there exists such a characteristic subgroup $H^*$ in $G$ that $B\cap H^*=H$ and $H$ is basic in $H^*$.

(2) If $H$ is an unbounded characteristic subgroup of a separable $p$-group $G$ and $B$ is basic in $G$, then $B\cap H$ is an unbounded characteristic subgroup in $B$ and $B\cap H$ is basic in $H$.
\end{lemma}

\Pf (1) Let $H^*=\sum_{f\in\mathrm{Aut}(G)}f(H)$ and $x\in B\cap H^*$. Then $x=y_1+\dots+y_n$, where $y_i=f_i(x_i)$,
$f_i\in\mathrm{Aut}(G)$, $x_i\in H$; $i=1, \dots, n$. Since $G$ is separable, the elements $x, y_1,\dots,y_n$ can be included in some basic subgroup $B'$ of $G$ (see, e.g., \cite{F1}). Let $\varphi: B\to B'$ be an isomorphism. However, because $G$ is torsion-complete, then $\varphi$ could be considered as an automorphism of $G$ (cf. \cite{F1}). Therefore, $z=\varphi^{-1}(x)$, $z_i=\varphi^{-1}(y_i)\in B$ and $U(x)=U(z)$, $U(z_i)=U(x_i)$. Since
$z_i=\psi_i(x_i)$ for some $\psi_i\in\mathrm{Aut}(B)$, we get that $z_i\in H$. So, $z=z_1+\dots+z_n\in H$. Furthermore, since $U(z)=U(x)$, it must be that $x\in H$. But Lemma~\ref{13} yields that $B\cap H^*=H$ is basic in
$H^*$, as expected.

(2) Again the application of Lemma~\ref{13} is a guarantor that we just need to show that $B\cap H$ is an unbounded characteristic subgroup in $B$. But, to that purpose, we observe that $B\cap H$ is really unbounded, because $H$ is unbounded and $B\cap H$ is basic in $H$. On the other hand, $B\cap H$ is characteristic in $B$ owing to Lemma~\ref{14}, as promised.
\fine

And so, we arrive at the next crucial result which, actually, is the converse of Theorem~\ref{15} in the case of torsion-complete groups. Specifically, the following assertion is true:

\begin{theorem}\label{17} A torsion-complete $p$-group $G$ is a weakly IC-group if, and only if, some (and hence all) of its basic subgroups $B$ is (are) a weakly IC-group.
\end{theorem}

\Pf Referring to Theorem~\ref{15}, one observes that it is necessary to prove only the necessity, because we already established the sufficiency. However, we shall demonstrate below slightly different arguments to show once again its validity.

To prove necessity, let $H$ be characteristic in $G$ with $H\neq G$. Since $G$ is a weakly IC-group and $H\cong G$, we deduce that $H$ as well as $G$ is unbounded. But then $H\cap B$ is characteristic in $B$ and $H\cap B$ is basic in $H$ in accordance with Lemma~\ref{16} (2). Note that $H\cap B\neq B$, because if we assume that $B\leq H$, then in view of the equality $B+H=G$ it would follow that $H=G$ - a contradiction. But since $H\cong G$, it follows that their basic subgroups have to be isomorphic (see, for instance, \cite{F}), say $B\cap H\cong B$, i.e., $B$ is also a weakly IC-group, as formulated.

To show now the truthfulness of sufficiency, let us assume that $G$ be a torsion-complete $p$-group and $B$ its basic subgroup which is a weakly IC-group. According to Proposition~\ref{03}, $B$ is an unbounded group and, therefore, $G$ is also unbounded. Since $B$ is a weakly IC-group, there exists a proper characteristic subgroup $S$ of the group $B$ such that $B\cong S$. Owing to Lemma~\ref{16}, there exists a proper characteristic subgroup $S^*$ of $G$ such that $S^*\cap B = S$ and $S$ is a basic subgroup of the group $S^*$. If the quotient $D/S^*$ is the divisible part of the quotient $G/S^*$, then a routine verification shows that $D$ is also a characteristic subgroup of $G$ and $S$ is a basic subgroup in $D$. Since the factor-group $G/D$ is reduced, it follows that $D$ is a torsion-complete group (see, e.g., \cite[Corrolary 68.7]{F}). Thus, we obtain that the group $G$ has a proper characteristic subgroup $D$ such that the basic subgroup $B$ of the group $G$ is isomorphic to the basic subgroup $S$ of the group $D$. Finally, as both $G$ and $D$ are torsion-complete groups, we conclude with the aid of volume II from \cite{F} that $G\cong D$, i.e., $G$ is a weakly IC-group, as required.
\fine

The following commentaries are worth noticing.

\medskip

The property "torsion-completeness" plays a key if not facilitating role in proving the last statement. However, a problem which automatically arises is of whether or not the theorem remains true by replacing the requirement "a torsion-complete $p$-group" by "a separable $p$-group", that is, is the statement of Theorem~\ref{15} reversible? A logical reason for asking this is that, in view of volume II of \cite{F}, we know that any separable $p$-group can be embedded as a pure and dense subgroup of a torsion-complete $p$-group. We, however, conjecture that the answer will be negative.

\section{Concluding Discussion and Open Problems}

We conclude this final section with some valuable comments on the obtained results.

And so, we firstly will comment the Abelian groups with an isomorphic proper characteristic subgroup, which groups were termed here as weakly IC-groups. Some of the analogies and differences with the classes of IC-groups and strongly IC-groups, both introduced in \cite{CD}, are as follows: (1) These three classes are definitely {\it not} closed under taking direct summands; (2) They have totally different structure as shown above.

These discrepancies help us to understand more completely the complicated situations in these three possibilities.

\medskip

We close our work with two questions of interest and importance.

\medskip

In regard to Examples~\ref{direct1} and \ref{direct2}, we pose the following left query.

\medskip

\noindent{\bf Problem 1.} Find necessary (and sufficient) conditions under which a direct summand of an IC-group is again an IC-group. Same holds and for strongly/weakly IC-groups, respectively.

\medskip

In connection with Corollary~\ref{06}, we state the following left question.

\medskip

\noindent{\bf Problem 2.} Do there exist (strongly, weakly) IC-groups which are not (fully) transitive?

\medskip

%\noindent{\bf Acknowledgements.} The authors would like to thank the unknown referee for his/her valuable comments and suggestions which led to an improvement of the article's shape.

\medskip

\noindent{\bf Funding:} The scientific work of Andrey R. Chekhlov was supported by the Ministry of Science and Higher Education of Russia (agreement No. 075-02-2022-884). The scientific work of Peter V. Danchev was partially supported by the Bulgarian National Science Fund under Grant KP-06 No 32/1 of December 07, 2019 and by the Junta de Andaluc\'ia, FQM 264.

\end{document}